\documentclass[10pt]{article}
\usepackage{float,amsmath,amsthm,amssymb,
	latexsym,amsfonts,amstext,amscd,
	amsbsy,bezier,enumerate,bigints}
\theoremstyle{plain}

\newtheorem{Theorem}{Theorem}
\newtheorem{Remark}{Remark}
\newtheorem{Corollary}{Corollary}
\newtheorem{Example}{Example}
\newtheorem{Lemma}{Lemma}
\newtheorem{Definition}{Definition}
\newcommand{\cqd}{\hfill \rule{2mm}{2mm}}

\title{{\bf\Large Delta-epsilon functions and uniform continuity on metric spaces}}
\author{{\bf\large C\'esar A. Hern\'andez M.}\footnote{Email: cahmelo@uem.br, Tel. (55)-44-3011-5358, Fax: (55)-11-3011-3873 }\hspace{2mm}
	{\bf\large}\vspace{1mm}\\
	{\it\small Department of Mathematics,  DMA-UEM}\\
	{\it\small Av. Colombo, 5790 Jd. Universit\'ario,,}\\
	{\it\small  CEP 87020-900, Maring\'a, PR, Brazil}}
	\vspace{3mm}
\date{October 9, 2017}
\begin{document}
	\maketitle
	\begin{abstract}
	Under certain general conditions, an explicit formula to compute the greatest delta-epsilon function of a continuous function is given. From
this formula, a new way to analyze the uniform continuity of a continuous function is given. Several examples illustrating the theory are discussed. 
\end{abstract}
	\textbf{Key  words}. Delta-epsilon functions, Continuity, Uniform continuity.

\section{Introduction}\label{epsilon}

Directly or indirectly most results in mathematical analysis use the concept of continuity in order to extend a property of a function $f$ that is satisfied at a point $p$ to a property satisfied in a neighborhood of $p$. A well known example that illustrates that fact is the Inverse function theorem, which we recall here:
\begin{itemize}
\item[] Roughly speaking, the Inverse function theorem states that a continuously differentiable mapping $f:\mathbb{R}^n\rightarrow \mathbb{R}^n$ is invertible on a certain open ball $B(p,\delta)$ where the linear transformation $f'(p)$ is invertible.  
\end{itemize}

In this example, the invertibility of $f'$ at $p$ is extended to an open ball via the continuity of the mapping $x\rightarrow f'(x)$ at $p$. Now, we notice that in the statement of the Inverse function theorem, nothing is said about the size of the ball where the function $f$ is invertible.  
However, it is not difficult to check, see for instance \cite{Rudin}, that the radius of the open ball depends on the norm of the linear transformation $[f'(p)]^{-1}$ and on the positive number delta appearing in the definition of continuity of the function $x\rightarrow f'(x)$ at the point $p$. More exactly, for $2\lambda||[f'(p)]^{-1}||=1$, then $\delta$ is such that 
$$
\text{If,}\hspace{0.5cm}||x-p||<\delta,\hspace{0.5cm}\text{then},\hspace{0.5cm}||f'(x)-f'(p)||<\lambda.
$$
So, inspired by the previous discussion, we would like to deal with the following questions:
\begin{enumerate}
\item[1.]\label{Otrica} Let $X,Y$ be metric spaces, $p\in X$ and $f:X\to Y$ a continuous function at $p$, and let $\epsilon>0$ be fixed, which is the greatest positive $\delta(p,\epsilon)$ such that the epsilon-delta definition of continuity is satisfied? Is there a formula to compute it?
\end{enumerate}

The previous questions has been solved indirectly in \cite{Hernandez} under simple assumptions on the function $f$. In fact, considering 
\begin{equation}\label{efe}
f:[b,\infty)\rightarrow [f(b),\infty)\hspace{0.4cm}\text{an increasing bijective function,}
\end{equation} 
$p\in[b,\infty)$, $\epsilon>0$ fixed, $p_0=f^{-1}(f(b)+\epsilon)$, $\delta_1(p,\epsilon)=p-f^{-1}(f(p)-\epsilon)$, and $\delta_2(p,\epsilon)=f^{-1}(f(p)+\epsilon)-p$, then the positive number  
\begin{equation}\label{deltinha}
\delta(p,\epsilon)=
\begin{cases}
\delta_2(p,\epsilon) & b\leq p<p_0\\
\min\{\delta_1(p,\epsilon),\delta_2(p,\epsilon)\} & p_0\leq p
\end{cases}
\end{equation}
satisfies the definition of continuity of $f$ at the point $p$. More precisely
\begin{equation}\label{conf}
\text{If,}\hspace{0.3cm}x\in[b,\infty)\hspace{0.2cm}\text{and}\hspace{0.2cm} |x-p|<\delta(p,\epsilon),\hspace{0.3cm}\text{then}\hspace{0.3cm}|f(x)-f(p)|<\epsilon.
\end{equation}In addition, $\delta(p,\epsilon)$ is the maximum positive number satisfying the previous condition. One of the aims of this paper is to generalize the formula (\ref{deltinha}) to the case of functions defined on metric spaces.\\

Now, it is also shown in \cite{Hernandez}, that the function $p\in[b,\infty)\rightarrow\delta(p,\epsilon)$ in (\ref{deltinha}) provides a way to study the uniform continuity of increasing bijective functions defined on unbounded intervals. Indeed, it is proven in \cite{Hernandez} that an increasing bijective function $
f:[b,\infty)\rightarrow [f(b),\infty)$ is not uniformly continuous on $[b,\infty)$, 
if and only if, 
there exists $\epsilon_0>0$, such that 
$$
\displaystyle{\inf_{x\in[b,\infty)}}\delta(x,\epsilon_0)=0.
$$
\\
Thus, other natural questions we are interested in are the followings: 
\begin{enumerate}
\item[2.]\label{Otra} Is it possible to study the uniform continuity of a continuous function $f:X\rightarrow Y$ in terms of the function $p\in X\rightarrow\delta(p,\epsilon)$? Are there other mathematical problems where the function $p\in X\rightarrow\delta(p,\epsilon)$ is useful? 
\end{enumerate}

Finally, it is important to point out that the question 1 and 2 are related to the concept of modulus of continuity
which has been extensively used in approximation theory, 
see for instace \cite{GSG}, \cite{AP} and the references therein.\\ 

This manuscript will be divided as follows, in the section \ref{reta}, we deal with the question \ref{Otrica}, more exactly, for continuous functions defined on certain metric spaces, an explicit formula to compute the greatest positive number $\delta(p,\epsilon)$ is given, see theorem \ref{Delta} and its corollaries. On the other hand, in the section \ref{uniform}, the question \ref{Otra} is discussed. Specifically, the theorems 2.1 and 2.2 proved in \cite{Hernandez} are generalized, see theorems \ref{intuitive} and \ref{beleza}. Additionally, some examples are discussed in the section \ref{Four}. 

\section{A formula to compute delta-epsilon numbers}\label{reta}

In this section, we give an explicit formula to compute delta-epsilon numbers of continuous functions. 
\begin{Definition}\label{definition}
Let $X,Y$ be metric spaces, $f:X\rightarrow Y$ a continuous function at $p\in X$, and $\epsilon>0$. A positive number $\delta$ is said to be a \textbf{delta-epsilon number} for $f$ at $p$, if $\delta$ satisfies the $\epsilon$-$\delta$ definition of continuity of $f$ at the point $p$. In other words, $\delta$ is such that
\begin{equation}\label{def}
\text{if}\hspace{0.2cm} x\in X\hspace{0.2cm}\text{and}\hspace{0.2cm}d_X(x,p)<\delta, \hspace{0.4cm}\text{then}\hspace{0.4cm}d_Y(f(x),f(p))<\epsilon.
\end{equation}
\end{Definition}
The following theorem provides a theoretical formula to compute the greatest delta-epsilon number for a wide class of functions. Namely, 
\begin{Theorem}\label{Delta}Let $f:X\rightarrow Y$ be a continuous function on $X$, $p\in X$ and $\epsilon>0$, then we have:
\begin{enumerate}
\item\label{Item1} If $f^{-1}(S[f(p),\epsilon])\neq\emptyset$, then the quantity 
\begin{equation}\label{Formula1}
\delta(p,\epsilon)=\text{dist}(p,f^{-1}(S[f(p),\epsilon])),
\end{equation}
is well defined and represents a positive number. Here, $S[f(p),\epsilon]$ denotes the sphere with center at $f(p)$ and radius $\epsilon$, thats is to say, $S[f(p),\epsilon]=\{y\in Y|d_Y(f(p),y)=\epsilon\}$.
\item\label{Item2}Furthermore, if the open ball $B(p,\delta(p,\epsilon))$ is path-connected then the number $\delta(p,\epsilon)$ is a delta-epsilon number for $f$ at $p$. More exactly, for every $x\in X$ such that, $d_X(x,p)<\delta(p,\epsilon)$, then $d_Y(f(x),f(p))<\epsilon$. 
\item\label{Item3} $\delta(p,\epsilon)$ is the greatest delta-epsilon number at $p$.
\item\label{Item4} Finally, if we define the set $D_{p,\epsilon}$ as:
\begin{equation}\label{D}
D_{p,\epsilon}=\{\beta\in\mathbb{R}^+|(\forall x\in X)(d_X(x,p)<\beta\Rightarrow d_Y(f(x),f(p))<\epsilon)\},
\end{equation}
then, $\delta(p,\epsilon)=\max D_{p,\epsilon}$ and of course $D_{p,\epsilon}=(0,\delta(p,\epsilon)]$.
\end{enumerate}
\end{Theorem}
\begin{Proof} The proof of this theorem proceed as follows:
\begin{enumerate} 
\item We first observe that since $f^{-1}(S[f(p),\epsilon])$ is a nonempty set, then the number  
\begin{equation}\label{Formula2}
\delta(p,\epsilon)=\inf\{\text{ }d_X(x,p)\text{ }| x\in X,\text{ }d_Y(f(x),f(p))=\epsilon\text{ }\},
\end{equation}
is well defined. Now, if $\delta(p,\epsilon)=0$, there exists a sequence $x_n\in X$ so that $\lim d_X(x_n,p)=0$ with $\lim d_Y(f(x_n),f(p))=\epsilon$. Being $f$ continuous at $p$, we can conclude that $\lim d_Y(f(x_n),f(p))=0$, since $\epsilon>0$, we have a contradiction. Thus, $\delta(p,\epsilon)$ have to be a positive number.
\item Let $x\in X$, we want to prove that,
\begin{equation}
\text{if},\hspace{0.3cm}d_X(x,p)<\delta(p,\epsilon)\hspace{0.3cm}\text{then}\hspace{0.3cm}d_Y(f(x),f(p))<\epsilon.
\end{equation} 
In fact, because of the definition of $\delta(p,\epsilon)$, clearly $d_Y(f(x),f(p))\neq\epsilon$. So, to finish the proof of our statement, we must show that the inequality $d_Y(f(x),f(p))>\epsilon$ is not possible. Now, we argue by contradiction. If $d_Y(f(x),f(p))>\epsilon$, since the open ball $B(p,\delta(p,\epsilon))$ is path-connected, there exists a continuous function $\gamma:[0,1]\rightarrow B(p,\delta(p,\epsilon))$ such that $\gamma(0)=p$ and $\gamma(1)=x$. Therefore, the function $g:[0,1]\rightarrow\mathbb{R}$ given by $g(t)=d_Y(f(\gamma(t)),f(p))$ is continuous and satisfies that $g(0)=0$ and $g(1)>\epsilon$, so by the intermediate value theorem, there exists $t_0\in(0,1)$ such that $g(t_0)=d_Y(f(\gamma(t_0)),f(p))=\epsilon$. Thus, $\gamma(t_0)$ satisfies that $d_X(\gamma(t_0),p)<\delta(p,\epsilon)$ and $d_Y(f(\gamma(t_0)),f(p))=\epsilon$, the last affirmation contradicts the definition of $\delta(p,\epsilon)$. Hence, we can conclude that the number $\delta(p,\epsilon)$ given in (\ref{Formula1}) is a delta-epsilon number for $f$ at $p$. As we wanted to prove.
\item If $\alpha$ is such that $\delta(p,\epsilon)<\alpha$ then by definition of $\delta(p,\epsilon)$ there exists $x\in X$ so that $\delta(p,\epsilon)\leq d_X(x,p)<\alpha$ with $ d_Y(f(x),f(p))=\epsilon$, so $\alpha$ is not a delta-epsilon number for $f$ at $p$. 
\item Now we proceed to prove item \ref{Item4}. From items \ref{Item1} and \ref{Item2}, we deduce that $\delta(p,\epsilon)\in D_{p,\epsilon}$. From item \ref{Item3}, we obtain that any other number greater than $\delta(p,\epsilon)$ is not in $D_{p,\epsilon}$. Hence, we can conclude that $\delta(p,\epsilon)=\max D_{p,\epsilon}$. This finishes the proof of the theorem.
\end{enumerate}
\cqd
\end{Proof}\\

\begin{Remark}
 In the item \ref{Item2} of the previous theorem, the connexity of the open ball $B(p,\delta(p,\epsilon))$ is a sufficient condition to get that $\delta(p,\epsilon)$ is a delta-epsilon number for $f$ at $p$. In fact, if there exists $x_0\in B(p,\delta(p,\epsilon))$ such that $d_Y(f(x_0),f(p))>\epsilon$, then, the  function $h$ given by
 $$
 h:B(p,\delta(p,\epsilon))\rightarrow \mathbb{R},\hspace{0.5cm} h(x)=d_Y(f(x),f(p))
 $$  
is continuous and satisfies $h(p)=0<\epsilon<h(x_0)$, then by the intermetiate value theorem, there exists $x_1\in B(p,\delta(p,\epsilon))$ such that $d_Y(f(x),f(p))=\epsilon$, which contradicts the definition of $\delta(p,\epsilon)$. 
\end{Remark}

The following results give us sufficient conditions to compute delta-epsilon numbers for points $p,x\in X$ which are connected by a path and such that $f(p)\neq f(x)$. More accurately,
\begin{Corollary}\label{Corollary1}
Let $f:X\rightarrow Y$ be a continuous function on $X$, suppose that there exist
$p,x\in X$ such that $d_Y(f(x),f(p))=:\beta>0$ and there exists a path connecting the points $p$ and $x$, then for every $\epsilon$ such that $0<\epsilon\leq\beta$ we have that $f^{-1}(S[f(p),\epsilon])\neq\emptyset$ and $f^{-1}(S[f(x),\epsilon])\neq\emptyset$. Particularly, for every $\epsilon$ satisfying $0<\epsilon\leq\beta$, the numbers $\delta(p,\epsilon)$ and $\delta(x,\epsilon)$ given by the formula (\ref{Formula1}) are well defined and positives. 
\end{Corollary}
\begin{Proof}
Since there exists a path $\gamma:[0,1]\rightarrow X$ connecting $p$ and $x$, then the function $g(t):[0,1]\rightarrow \mathbb{R}$ given by $d_Y(f(\gamma(t),f(p)))$ is continuous   
and satisfies $[0,\beta]\subset g([0,1])$. Thus, for every $\epsilon$ such that $0<\epsilon<\beta$ there exists $t_0\in(0,1)$ satisfying $d_Y(f(\gamma(t_0),f(p)))=\epsilon$, which proves that $f^{-1}(S[f(p),\epsilon])\neq\emptyset$. Similarly, it is shown that $f^{-1}(S[f(x),\epsilon])\neq\emptyset$. The rest of the proof follows from the item \ref{Item1} in the theorem \ref{Delta}.
\cqd
\end{Proof}
\begin{Corollary}\label{Corollary2}
Let $f:X\rightarrow Y$ be a continuous function on $X$, suppose that there exist
$p,x\in X$ such that $d_Y(f(x),f(p))=:\beta>0$, and there exists a path connecting the points $p$ and $x$, if for every $\epsilon$ with $0<\epsilon\leq\beta$, the open balls $B(p,\delta(p,\epsilon))$, $B(x,\delta(x,\epsilon))$  are path-connected, then the numbers $\delta(p,\epsilon)$, $\delta(x,\epsilon)$ are delta-epsilon numbers for $f$ at $p$ and $x$ respectively.  
\end{Corollary}
\begin{Proof}
The proof of this result follows from the corollary \ref{Corollary1} and the proof of the item \ref{Item2} in the theorem \ref{Delta}. 
\cqd
\end{Proof}\\

The following theorem allows us to compute delta-epsilon numbers in a neighborhood of a point $p$ which admits a delta-epsilon number. Namely,
\begin{Theorem}\label{Open}
Let $f:X\rightarrow Y$ be a continuous function on $X$, suppose that there exist
$p,x\in X$ such that $d_Y(f(x),f(p))=:\beta>0$, suppose that the open ball $B(p,\delta(p,\beta))$ is path connected and that there exists a path connecting the points $p$ and $x$.
Then, for every $\epsilon$, with $0<\epsilon<\beta$, there exists $\delta$ satisfying $0<\delta\leq\delta(p,\beta)$, such that if $d_X(q,p)<\delta$ then the numbers $\delta(q,\epsilon)$ given in (\ref{Formula1}) are well defined an positives. If in addition, the open balls $B(q,\delta(q,\epsilon))$ are path-connected, then for all $q\in B(p,\delta)$, the numbers $\delta(q,\epsilon)$ are delta-epsilon numbers.  
\end{Theorem}
\begin{Proof} We divide the proof of this theorem into two parts:

\begin{enumerate}
\item First all of, we shall prove that there exists $\delta$ with $0<\delta\leq\delta(p,\beta)$ and such that if $d_X(q,p)<\delta$ then $\epsilon<d_Y(f(x),f(q))$. In fact, since $f^{-1}(S[f(p),\beta])\neq\emptyset$ and the open ball $B(p,\delta(p,\beta))$ is path connected, then from the theorem \ref{Delta}, we conclude that the number $\delta(p,\beta)$ is the maximum delta-epsilon number at $p$. On the other hand, since $f$ is continuous at $p$ and $\beta-\epsilon$ is positive, there exists $\delta>0$ such that if $d_X(q,p)<\delta$ then $d_Y(f(q),f(p))<\beta-\epsilon<\beta$. So, since $\delta(p,\beta)$ is the maximum delta-epsilon number at $p$, we deduce that $\delta\leq\delta(p,\beta)$. Now, by taking $q\in B(p,\delta)$ and from triangular inequality, we obtain that 
\begin{equation}
\begin{aligned}
\beta=d_Y(f(x),f(p))&\leq d_Y(f(x),f(q))+d_Y(f(q),f(p))
\\&<d_Y(f(x),f(q))+\beta-\epsilon,
\end{aligned}
\end{equation}
so, if $d_X(q,p)<\delta$ then $\epsilon<d_Y(f(x),f(q))$. As we wanted to show. 
\item Final part, the conclusion of the proof. As each point $q$ in the ball $B(p,\delta)$ can be connected with $x$ using a path and since $\epsilon<d_Y(f(x),f(q))$, then the corollary \ref{Corollary1} leads us to conclude that $f^{-1}(S[f(q),\epsilon])\neq\emptyset$. So the numbers $\delta(q,\epsilon)$ are well defined in the ball $B(p,\delta)$. Finally, since the open balls $B(q,\delta(q,\epsilon))$ are path-connected then the item \ref{Item2} in the theorem \ref{Delta} allows us to conclude that the numbers $\delta(q,\epsilon)$ are delta-epsilon numbers. This finishes the proof of the theorem. 
\end{enumerate}
\cqd
\end{Proof}\\

Finally, the next corollary gives us sufficient conditions to calculate delta-epsilon numbers in all of the domain of the function $f$.
\begin{Corollary}\label{Corollary3}
Let  $f:X\rightarrow Y$ be a nonconstant continuous function defined on a metric space $X$. If for all $p\in X$ and $r>0$ the open balls $B(p,r)$ are path-connected, then there exists $\beta>0$ such that, the delta-epsilon numbers $\delta(p,\epsilon)$ are well defined on the set $X\times(0,\beta)$. 
\end{Corollary}
\begin{Proof} 
According to the corollary \ref{Corollary1}, to show this result, it is necessary to find out a positive number $\beta$ such that for every $p\in X$ there exists $x\in X$ satisfying that $d_Y(f(p),f(x))=\beta$. In fact, since $f$ is a nonconstant function, then the diameter of $f(X)$ is positive, namely, $diam(f(X))>R$ for some $R>0$. So, there exist $a,b\in X$ with $R/2<d_Y(f(a),f(b))$. Now, for $p\in X$, we have that $$R/2<d_Y(f(a),f(b))\leq d_Y(f(a),f(p))+d_Y(f(p),f(b)),$$ thus, we can conclude that either $R/4<d_Y(f(a),f(p))$ or $R/4<d_Y(f(p),f(b))$. On the other hand, since $X$ is path-connecting, then there exists $x\in X$ such that $d_Y(f(x),f(p))=R/4$. Finally, the proof of the corollary follows from direct application of the corollaries \ref{Corollary1}, \ref{Corollary2} by taking $\beta:=R/4$.   
\cqd
\end{Proof}\\

Now, we establish some properties of the delta-epsilon numbers,
\begin{Theorem}\label{Properties}
Let  $f:X\rightarrow Y$ be a nonconstant continuous function defined on a compact metric space $X$. Suppose that for all $p\in X$, $r>0$ the open balls $B(p,r)$ are path-connected,  let $\beta$ be the positive number obtained in the corollary \ref{Corollary3}, then the function $\delta:X\times(0,\beta)\rightarrow \mathbb{R}^+$ defined by $(p,\epsilon)\rightarrow \delta(p,\epsilon)$ satisfies the following properties:
\begin{enumerate}
\item\label{A} For all $p\in X$ and $\epsilon>0$ there exists $x\in X$ such that $\delta(p,\epsilon)=d_X(p,x)$, and $d_Y(f(p),f(x))=\epsilon$.
\item\label{B} If $a<b<\beta$, then $\delta(p,a)\leq\delta(p,b)$.
\item\label{D} $\lim_{n\to\infty} \delta(p,\epsilon-1/n)=\delta(p,\epsilon)$.
\item\label{C} Let $x_n$ a sequence in $X$ such that $\lim_{n\to\infty} x_n=p$, then for all $r>0$, there exists $n_0\in\mathbb{N}$, so that for all $n>n_0$ $$\delta(p,\epsilon-r)-r\leq\delta(x_n,\epsilon)\leq\delta(p,\epsilon+r)+r.$$
\end{enumerate}
\end{Theorem}
\begin{Proof}
\begin{enumerate}
\item Since $f$ is a nonconstant continuous function, and $X$ is a compact path-connected set, then $f^{-1}(S[f(p),\epsilon])$ is a nonempty compact set that not contains $p$. So the distance between $p$ and $f^{-1}(S[f(p),\epsilon])$ is reached at some point $x\in f^{-1}(S[f(p),\epsilon])$. That is to say,
$$\delta(p,\epsilon)=dist(p,f^{-1}(S[f(p),\epsilon])=d_X(p,x),$$ 
with $d_Y(f(p),f(x))=\epsilon$.
\item In terms of the notation of the theorem \ref{Delta}, since $a<b$, then $\delta(p,a)\in D_{p,b}$, and since $\delta(p,b)=\max D_{p,b}$, then we conclude that $\delta(p,a)\leq\delta(p,b)$.
 \item From the item \ref{B}, the sequence $\delta(p,\epsilon-1/n)$ is increasing and bounded from above by $\delta(p,\epsilon)$. On the other hand, from the item \ref{A}, there exists a sequence $x_n\in X$ such that 
$$
\delta(p,\epsilon-1/n)=d_X(p,x_n),\hspace{0.3cm}\text{and}\hspace{0.3cm}d_Y(f(p),f(x_n))=\epsilon -\frac{1}{n}.
$$
Now, since $X$ is compact, then there exists $q\in X$ such that $\lim_{k\to\infty}x_{n_k}=q$ where $x_{n_k}$ is a subsequence of the sequence $x_n$. So, by the continuity of the function $f$, we have that $\lim_{n\to\infty}\delta(p,\epsilon-1/n)=d_X(p,q)$ with $d_Y(f(p),f(q))=\epsilon$, then we have tha $d_X(p,q)\leq\delta(p,\epsilon)\leq d_X(p,q)$, thus, we obtain that $\lim_{n\to\infty} \delta(p,\epsilon-1/n)=\delta(p,\epsilon)$. This finishes the proof of item \ref{D}.
\item Since $f$ is continuous at $p\in X$ and $\lim_{n\to\infty}x_n=p$, then $\lim_{n\to\infty}f(x_n)=f(p)$. So, for $r>0$ there exists $n_0\in\mathbb{N}$ such that for all $n>n_0$, $d_X(x_n,p)<r$ and $d_Y(f(x_n),f(p))<r$. Now, we prove that $\delta(p,\epsilon-r)-r\in D_{x_n,\epsilon}$. In fact, let $y\in X$ with $d_X(y,x_n)<\delta(p,\epsilon-r)-r$, then from triangular inequality, we obtain that $d_X(y,p)<\delta(p,\epsilon-r)$ and so, $d_Y(f(y),f(p))<\epsilon-r$. Finally, $d_Y(f(y),f(x_n))\leq d_Y(f(y),f(p))+d_Y(f(p),f(x_n))<\epsilon$. This finishes the proof of the first inequality of item \ref{C}. The proof of the second inequality can be done similarly.     

\end{enumerate}
\end{Proof}


\section{Uniform continuity and $\delta$-$\epsilon$ functions}\label{uniform} 

In this section, we extend the concept of \textit{delta-epsilon function} introduced in \cite{Hernandez} and use it to study the relationship between continuity and uniform continuity.  
\begin{Definition}\label{contunif}
Let $X$ be a nonempty set, a function $f:X\rightarrow Y$ is called uniformly continuous on $X$, if for every $\epsilon>0$, there exists $\delta>0$ such that 
for every $x,y\in X$ with $d_X(x,y)<\delta$, then $d_Y(f(x),f(y))<\epsilon$.  
\end{Definition}
\begin{Definition}\label{Delta-Epsilon}
Let $X$ be a nonempty set, $f:X\rightarrow Y$ a continuous function. Let $\epsilon>0$ fixed, we say that a function $g_{\epsilon}:X\rightarrow \mathbb{R}^+$ is a \textbf{delta-epsilon function} for $f$,  if $g_{\epsilon}(x)$ is a delta-epsilon number for $f$ at $x\in X$. 
\end{Definition}
\begin{Example}\label{Unico}\rm{
Let  $f:X\rightarrow Y$ be a nonconstant continuous function defined on a metric space $X$, and suppose that for all $p\in X$ and $r>0$ the open balls $B(p,r)\subset X$ are path-connected. Then from the corollary \ref{Corollary3} and the item \ref{Item3} in the theorem \ref{Delta}, we have that for any $\epsilon\in(0,\beta)$ fixed, the function   
$\delta(\cdot,\epsilon):X\rightarrow \mathbb{R}^+$, given by
\begin{equation}\label{hatem}
\delta(x,\epsilon)=\text{dist}(x,f^{-1}(S[f(x),\epsilon])), 
\end{equation} 
is a delta-epsilon function for $f$ that is greater that any other delta-epsilon function for $f$.
}
\end{Example}
The following theorem gives a characterization of the uniform continuity concept in terms of delta-epsilon functions,
\begin{Theorem}\label{intuitive} Let $X$ a nonempty set and $f:X\rightarrow Y$ a continuous function. Then, $f$ is continuous uniformly on $X$ if and only if there exists a family $\left\{g_{\epsilon}\right\}_{\epsilon>0}$ of delta-epsilon functions for $f$ such that, 
\begin{equation}\label{posi}
\eta_{\epsilon}:=\inf_{x\in X} g_{\epsilon}(x)>0,
\end{equation} for every $\epsilon>0$.  
\end{Theorem}
\begin{Proof} 
If $f:X\rightarrow Y$ is continuous uniformly on $X$, then for $\epsilon>0$ there exists $\delta>0$, such that 
for every $x,y\in X$ with $d_X(x,y)<\delta$, then  $d_Y(f(x),f(y))<\epsilon$. Thus, 
the constant function $g_{\epsilon}:X\rightarrow \mathbb{R}^+$, $g_{\epsilon}(x)=\delta$, is a delta-epsilon function for $f$ that clearly satisfies the condition (\ref{posi}).
\\Conversely, let $\left\{g_{\epsilon}\right\}_{\epsilon>0}$ a family of delta-epsilon functions for the continuous function $f$ that satisfies the condition (\ref{posi}). Hence, for $\epsilon>0$ and $x,y\in X$, we have that if $d_X(x,y)<\eta_{\epsilon}\leq g_{\epsilon}(x)$, then, since $f$ is continuous at $x$ and $g_{\epsilon}(x)$ satisfies the continuity definition at $x$, we can conclude that $d_Y(f(x),f(y))<\epsilon$.   
\cqd
\end{Proof}
\begin{Remark}\rm{
Roughly speaking, the previous theorem says us that a continuous function $f$ that admits a family of constant delta-epsilon functions $\left\{\eta_{\epsilon}\right\}_{\epsilon>0}$ is continuous uniformly.
}
\end{Remark}
\begin{Remark}\rm{
In terms of the theorem \ref{intuitive}, to show that a continuous function $f$ is not continuous uniformly, we must verify that any family $\left\{g_{\epsilon}\right\}_{\epsilon>0}$ of delta-epsilon functions has an element $g_{\epsilon_0}$  
such that, $\inf_{x\in X} g_{\epsilon_0}(x)=0$,
what seems to be a difficult task. For a certain class of functions the following theorem simplify this work.
}
\end{Remark}
Now, we give a characterization of the concept of uniform continuity in terms of the optimal  delta epsilon function given in (\ref{hatem}).
\begin{Theorem}\label{beleza}
Let $f:X\rightarrow Y$ be a nonconstant continuous function defined on a metric space $X$. Suppose that for all $p\in X$ and $r>0$ the open balls $B(p,r)$ are path-connected. Then, the following three conditions are equivalents:
\begin{enumerate}
\item\label{fc} $f$ is not uniformly continuous on X.
\item\label{co} There exists $\epsilon_0$ such that, $\displaystyle{\inf_{x\in X} \delta(x,\epsilon_0)=0}$.
\item\label{nv} There exist $\epsilon_0$ and sequences $x_n,y_n\in X$, such that, $\lim_{n\to\infty}d_X(x_n,y_n)=0$ and $d_Y(f(x_n),f(y_n))=\epsilon_0$. 
\end{enumerate}
\end{Theorem}
\begin{Proof}  
\begin{itemize}
\item First, we prove that \textit{\ref{fc}} implies \textit{\ref{co}}. It is clear that if $f$ is not continuous uniformly on $X$, then from the theorem \ref{intuitive}, the family of delta-epsilon functions $\left\{\delta(\cdot,\epsilon)\right\}_{\epsilon\in(0,\beta)}$ must have an element satisfying the condition \textit{\ref{co}}.
\item Now, we prove that \textit{\ref{co}} implies \textit{\ref{fc}}. Let $\left\{ \rho_{\epsilon}\right\}_{\epsilon>0}$ be a family of delta-epsilon functions for $f$.
Then, from the item \ref{Item3} of theorem \ref{Delta} and item \ref{B} of theorem \ref{Properties}, we have that
$$
\rho_{\epsilon}(x)\leq\delta(x,\epsilon)\leq \delta(x,\epsilon_0),\hspace{0.3cm}\text{fol all}\hspace{0.3cm} x\in X,
$$ 
where, $0<\epsilon\leq\epsilon_0$. Hence, from the condition \textit{\ref{co}} above, we obtain that $\inf_{x\in X}\rho_{\epsilon}(x)=0.$
So, from theorem \ref{intuitive}, we can conclude that $f$ is not continuous uniformly on $X$.
\item Next, we prove that \textit{\ref{co}} implies \textit{\ref{nv}}. Since $\inf_{x\in X}\delta(x,\epsilon_0)=0$, then for all $n\in\mathbb{N}$, there exists $x_n\in X$ such that 
$0<\delta(x_n,\epsilon_0)<1/n$. By the definition of $\delta(x_n,\epsilon_0)$ there exist $y_n\in X$, satisfying, $0<\delta(x_n,\epsilon_0)\leq d_X(x_n,y_n)<1/n$ and $d_Y(f(x_n),f(y_n))=\epsilon_0$. Thus, we obtain two sequences of elements $x_n,y_n\in X$, such that, $\lim_{n\to\infty}d_X(x_n,y_n)=0$ and $d_Y(f(x_n),f(y_n))=\epsilon_0$.
\item Finally, we prove that \textit{\ref{nv}} implies \textit{\ref{co}}. If the condition \textit{\ref{nv}} holds, then, we deduce that $0<\delta(x_n,\epsilon_0)\leq d_X(x_n,y_n)$. Hence, $\lim_{n\to\infty}\delta(x_n,\epsilon_0)=0$, what implies that $\inf_{x\in X}\delta(x,\epsilon_0)=0$. As we wanted to show. 
\cqd
\end{itemize}
\end{Proof}
\section{Examples}\label{Four}
In this section, from some specific examples, we illustrate the theory developed in sections \ref{reta} and \ref{uniform}. The first example is elementary, however, it explains our principal results clearly. 

\begin{Example}\rm{Consider $f:\mathbb{R}\rightarrow\mathbb{R}$ defined by $f(x)=x^2$, according to the theorem \ref{Delta}, it is clear that 
$$
f^{-1}(S[f(p),\epsilon])\supset\{\sqrt{p^2+\epsilon},-\sqrt{p^2+\epsilon}\}\neq\emptyset,
$$
then, the hypothesis 1 in the theorem \ref{Delta} is satisfied. After some calculations, we obtain that:
$$
\delta(p,\epsilon)=\text{dist}(p,f^{-1}(S[f(p),\epsilon]))=\sqrt{p^2+\epsilon}-|p|.
$$
Now, since $B(p,\delta(p,\epsilon))$ is connected, then, we can conclude that the numbers $\delta(p,\epsilon)$ are maximum delta-epsilon numbers for $f$,  that is:
\begin{equation*}
\text{if}\hspace{0.5cm}|x-p|<\sqrt{p^2+\epsilon}-|p|, \hspace{0.5cm}\text{then}\hspace{0.5cm} |x^2-p^2|<\epsilon. 
\end{equation*}
In addition, we obtain that
$$
D_{p,\epsilon}=\left(0,\sqrt{p^2+\epsilon}-|p|\right].
$$
On the other hand, as the function $p\to \delta(p,\epsilon)$ is even a decreassing for $p>0$, then, we deduce that for $M>0$
$$
\inf_{p\in B(0,M)}\delta(p,\epsilon)\geq\sqrt{M^2+\epsilon}-M>0 ,
$$
thus, from the theorem \ref{intuitive}, $f$ is continuous uniformly over bounded domains. In contrast, due to the fact that
$$
\lim_{|p|\to\infty}\delta(p,\epsilon)=0,
$$
then, from the theorem \ref{beleza}, we obtain that $f$ is not continuous uniformly over unbounded domains.
}
\end{Example}
The following lemmas allow us to study uniform continuity of radial functions defined on normed vector spaces,
\begin{Lemma}\label{interame}
Let $I\subset\mathbb{R}$ be a nondegenerate interval, $g:I\rightarrow\mathbb{R}$ an increasing bijective function. Then the numbers 
\begin{equation}\label{certa}
\delta(p,\epsilon)=\min\{g^{-1}(g(p)+\epsilon)-p, p-g^{-1}(g(p)-\epsilon)\}, 
\end{equation}
are maximum delta-epsilon numbers for the function $g$. Here, $\epsilon$ has to be chosen in such a way that some of the numbers $g^{-1}(g(p)+\epsilon)-p$, 
$p-g^{-1}(g(p)-\epsilon)$ be well defined.
\end{Lemma}
\begin{Proof}
Under these hypothesis, we can conclude that the function $g$ is continuous and nonconstant. Hence, since $I$ is connected, then the corollary \ref{Corollary3} implies the conclusion of the lemma. We point out that this lemma is exactly the lemma 1.2 in \cite{Hernandez}. 
\cqd
\end{Proof}

\begin{Lemma}\label{inter}
Let $V$ be a normed vector space, $I\subset\mathbb{R}^+$ be a nondegenerate interval, $g:I\rightarrow\mathbb{R}$ a nonconstant continuous function and $$\Omega=\{x\in V|\text{ } ||x||\in I\text{ }\},$$ if we define 
$f:\Omega\rightarrow\mathbb{R}$ by $f(x)=g(||x||)$ then, for $p\in\Omega$, we have that $\delta(||p||,\epsilon)\in D_{p,\epsilon}$. That is, the maximum delta-epsilon number $\delta(||p||,\epsilon)$ for $g$ at $||p||$ is a delta-epsilon number for the radial function $f$ at $p$. 
 \end{Lemma}
\begin{Proof}
First of all, we note that since $g$ is nonconstant and it is defined on an interval, then the corollary \ref{Corollary3} allows us to conclude  
that the numbers $\delta(||p||,\epsilon)$ are well defined in the sense of the formula (\ref{hatem}). Now, suppose that $$||x-p||<\delta(||p||,\epsilon),$$ then, $|||x||-||p|||<\delta(||p||,\epsilon)$, and thus, $$|f(x)-f(p)|=|g(||x||)-g(||p||)|<\epsilon,$$
as we wanted to proof.
\cqd
\end{Proof}
\begin{Example}\rm{
Let $V$ be a normed vector space, if we consider the function $f:V\rightarrow\mathbb{R}$, defined by 
$$
f(p)=e^{||p||},
$$  
then, from the example 6 in \cite{Hernandez} and the previous lemma, we can infer that the family 
\begin{equation}\label{dois}
\eta(p,\epsilon)=\ln(e^{||p||}+\epsilon)-||p||,
\end{equation}
is a family of delta-epsilon functions for $f$. In addition, since for $p\neq 0$, the function 
$$
t\to g(t)=e^{t||p||},
$$
is an homeomorphism from $[0,\infty)$ to $[1,\infty)$, then from the lemma \ref{interame}, we obtain that
$$
\delta(t,\epsilon)=\frac{\ln(e^{t||p||}+\epsilon)-t||p||}{||p||},
$$
is a family of maximum delta-epsilon functions for $g$. Now, as
$$
\lim_{t\to\infty} \frac{\ln(e^{t||p||}+\epsilon)-t||p||}{||p||}=0,
$$
we conclude that
$$
\inf_{t\in[0,\infty)} \delta(t,\epsilon)=0,
$$
and from the theorem \ref{beleza}, the function $g$ is not continuous uniformly on $[0,\infty)$. As a result of this analysis, we also obtain that the function $f$ is not continuous uniformly on any domain containing the set $\{tp|t\in[0,\infty)\}\subset V$. On the other hand, it is not difficult to show that the delta-epsilon functions given in (\ref{dois}) satisfy
$$
\inf_{p\in B(0,M)} \eta(p,\epsilon)\geq \ln(e^{M}+\epsilon)-M>0,
$$
hence, from the theorem \ref{intuitive}, we deduce that $f$ is continuous uniformly on bounded domains.
}
\end{Example}

\begin{Example}\rm{
Similarly, from the example 7 in \cite{Hernandez} and the lemma \ref{inter}, we can deduce that the function 
$f:V-\{0\}\rightarrow\mathbb{R}$, defined by 
$$
f(p)=\ln(||p||),
$$  
has the following family of delta-epsilon functions:
\begin{equation}\label{treis}
\eta(p,\epsilon)=||p||(1-e^{-\epsilon}),
\end{equation}
thus, from the theorem \ref{intuitive}, we obtain that for any $M>0$, $f$ is continuous uniformly on the domain $V-B(0,M)$. On the other hand, following the ideas of the previous example it is possible to conclude that $f$ is not continuous uniformly on $V-\{0\}$.  
}
\end{Example}

\renewcommand{\refname}{References}

\end{document}